\documentclass[reqno]{amsart}
\usepackage{amsmath, amsfonts,epsfig, amssymb}

{\theoremstyle{definition}
\newtheorem{definition}{Definition}[section]
\newtheorem{example}[definition]{Example}
}
{\theoremstyle{plain}%
  \newtheorem{theorem}{Theorem}
  \newtheorem{corollary}{Corollary}
  \newtheorem{proposition}[definition]{Proposition}
  \newtheorem{lemma}[definition]{Lemma}%
  \newtheorem{step}{Step}
}
{\theoremstyle{remark}
\newtheorem*{question}{Question}
\newtheorem{remark}[definition]{Remark}
}

\newcommand{\zz}{\mathbb{Z}}
\newcommand{\sn}{S_n}
\newcommand{\snelling}{snelling}
\newcommand{\snellable}{snellable}
\newcommand{\rk}{\mathrm{rk}}
\newcommand{\edges}{\mathcal{E}}
\newcommand{\ch}{\mathrm{ch}}
\newcommand{\mcp}{\mathcal{M}(P)}
\newcommand{\ghna}{good $\mathcal{H}_n(0)$ action}
\newcommand{\wm}{\omega_\m}
\newcommand{\shape}{bowtie}
\newcommand{\jpm}{J(P_{\omega_\m})}
\newcommand{\restless}{restless}
\newcommand{\m}{\mathfrak{m}}
\newcommand{\chain}{\mathfrak{c}}
\newcommand{\M}{\mathcal{M}}
\newcommand{\Mm}{\mathcal{M}_\m}

\newcommand{\subqed}{\quad\hbox{\framebox[2pt]{\rule{0pt}{3pt}}}}

\begin{document}

\title{EL-labelings, Supersolvability and 0-Hecke Algebra Actions on Posets}

\author{Peter McNamara}

\address{Department of Mathematics 2-342, Massachusetts Institute of
Technology, 77 Massachusetts Avenue, Cambridge, MA 02139, U.S.A.}

\email{mcnamara@math.mit.edu}


\begin{abstract}  
It is well known that if a finite graded lattice
of rank $n$ is supersolvable, 
then it has an EL-labeling where the labels along any maximal chain form 
a permutation.  We call such a labeling an $\sn$ EL-labeling and we show 
that a finite graded lattice of rank $n$ is supersolvable 
if and only if it has such a labeling.  We next consider finite graded posets 
of rank $n$ with $\hat{0}$ and $\hat{1}$ that have an $\sn$ 
EL-labeling.  We describe a type $A$ 0-Hecke 
algebra action on the maximal chains of such posets.  This action 
is local and gives a representation of these Hecke algebras whose character
has characteristic that is closely related to Ehrenborg's flag
quasisymmetric function.  We ask what other classes of posets have such 
an action and in particular we show that finite graded lattices of rank
$n$ have such an action if and only if they have an $\sn$ EL-labeling.
\end{abstract}



\maketitle


\section{Introduction}
Supersolvable lattices were introduced by R. Stanley in \cite{St2} where
he showed that the covering relations can be labeled by the 
integers to give an 
EL-labeling.  We explain and discuss these terms in
Section \ref{sec-snellable}.  In fact, this EL-labeling of a supersolvable
lattice of rank $n$ 
is seen to have the additional property that the labels along
any maximal chain of the lattice form a permutation of $1,2,\ldots ,n$.
We call this type
of labeling an $\sn$ EL-labeling.  In Section \ref{sec-proof1}, we prove 
that the converse result is true: if a finite graded lattice has 
an $\sn$ EL-labeling then it is supersolvable.

In Section \ref{sec-action}, we describe an action on the maximal chains
of an $\sn$ EL-labeled lattice, suggested to the author by R. Stanley.  
We show that this action
gives a representation of the Hecke algebra of type $A$ at $q=0$.  
In \cite{SS} and \cite{St3}, the Frobenius characteristic of the 
character of some symmetric group actions is shown to be closely
related to Ehrenborg's flag symmetric function.  In Section
\ref{sec-goodaction}, we show that our $\mathcal{H}_n(0)$ action
has an analogous property and we follow Simion and Stanley
in calling our action a \emph{good} $\mathcal{H}_n(0)$ action.
Note that the material of Section 
\ref{sec-goodaction} is not necessary for the  
proof of Section 
\ref{sec-proof1}. 
Our second main result appears in Section \ref{sec-proof2}.
We show that a certain class of posets, which includes finite
graded lattices, have a \ghna\ if and only if they have an $\sn$ EL-labeling.
It follows that a finite graded lattice is supersolvable if and only if
it has a \ghna.

\section{EL-labelings and Supersolvability}\label{sec-snellable}

Throughout, we let $s_i$ denote the permutation which transposes
$i$ and $i+1$, and composition of permutations will be from right to left.
For any positive integer $n$, write $[n]$ for the set $\{1, 2, \ldots, n\}$.
Suppose $P$ is a finite graded poset of rank $n$, 
with $\hat{0}$ and $\hat{1}$.
(For undefined poset terminology, see 
\cite[Ch.\ 3]{ec1}.)
Let $\rk$ denote the rank function of $P$, so 
$\rk(\hat{0})=0$ and 
$\rk(\hat{1})=n$.  
If $x \leq y$ in $P$, let $\rk(x,y)$ denote 
$\rk(y)-\rk(x).$ 
If $x \leq y$ in $P$ and $\rk(x,y)=1$ then we 
say that $y$ \emph{covers} $x$.  Let 
$\edges(P) = \{(s,t):t \mbox{ covers } s 
\mbox{ in } P\}$, the set of edges of the Hasse diagram of $P$, and let 
$\mcp$ denote the set of maximal chains of $P$. 

A function $\lambda : \edges(P) \to \zz$ gives us an \emph{edge-labeling}
of $P$.  If 
$\m : s = s_0 < s_1 < \cdots < s_k = t$
is a maximal chain
of the interval $[s,t]$, then we write 
$\lambda(\m) =
(\lambda(s_0, s_1), \lambda(s_1, s_2), \ldots , \lambda(s_{k-1}, s_k) ).$
The chain $\m$ is \emph{increasing} if 
\linebreak
$\lambda(s_0, s_1) \leq \lambda(s_1, s_2)
 \leq \cdots \leq \lambda(s_{k-1}, s_k).$ 
We let 
$\leq_{L}$ denote lexicographic 
order on finite integer 
sequences: 
$(a_1, a_2, \ldots ,a_k) <_{L} (b_1, b_2, \ldots ,b_k)$ 
if and only if
$a_i < b_i$ in the first coordinate where they differ.

\begin{definition}
Let $P$ be a finite graded poset of rank $n$.  An edge-labeling 
$\lambda : \edges(P) \to \zz$ is called an \emph{EL-labeling} if the
following two conditions are satisfied:
\begin{itemize}
\item[(i)] Every interval $[s,t]$ has exactly one increasing maximal chain $\m$.
\item[(ii)] Any other maximal chain $\m^{\prime}$ of $[s,t]$ satisfies
	$\lambda(\m^{\prime}) >_L \lambda(\m)$.
\end{itemize}
\end{definition}

A poset $P$ with an EL-labeling is said to be 
\emph{edge-wise lexicographically shellable} or \emph{EL-shellable}. 
This definition of a lexicographically shellable poset 
first appeared in \cite{Bj} with the motivating examples being
from \cite{St1} and \cite{St2}, which appear as Examples \ref{distributive}
and \ref{ss} below.  
The ubiquity and usefulness of EL-labelings arises from the fact that if
$P$ is EL-shellable, then $P$ is shellable and hence Cohen-Macaulay.
Further information on these concepts can be found in \cite{Bj} and 
the highly recommended survey article \cite{BGS}.
We will be interested in the following type of EL-labeling:

\begin{definition}
An EL-labeling $\lambda$ of $P$ is said to be an \emph{ $\sn$ EL-labeling }
if, for every maximal chain 
$\m : \hat{0}=x_0 < x_1 < \cdots < x_n=\hat{1}$ 
of $P$, the map sending $i$ to $\lambda(x_{i-1}, x_i)$ is a 
permutation of $[n]$.  In other words, $\lambda(\m)$ is a permutation of
$[n]$ written in the usual way. 
\end{definition}

If a poset $P$ has an $\sn$ EL-labeling, or \emph{\snelling} for short, 
then it is said to be \emph{$\sn$ EL-shellable}, or 
\emph{\snellable} for short.
Note that the second condition in the definition of an
EL-labeling is redundant in 
this case.  

\begin{example}\label{bn}
Consider the poset $B_n$, the set of subsets of $[n]$.  If $y$ covers $x$
in $B_n$ then $y - x=\{i\}$ for some $i \in [n]$ and we set $\lambda(x,y)= i$. 
This defines a \snelling\ for $B_n$.
\end{example}

\begin{example} \label{distributive}
Any finite distributive lattice is \snellable.
Let $L$ be a finite distributive lattice of rank $n$.  By 
\cite[p. 59, Thm. 3]{Bi}, that is 
equivalent to saying that
$L = J(Q)$, the lattice of order ideals of some $n$-element poset $Q$. 
Let $\omega : Q \to [n]$ be a linear extension of $Q$, i.e., any 
bijection labeling the vertices of Q that is
order-preserving (if $a<b$ in $Q$ then $\omega(a)<\omega(b)$). 
This labeling of the vertices of $Q$ defines a labeling of the edges of
$J(Q)$ as follows.  If $y$ covers $x$ in $J(Q)$, then the order ideal 
corresponding to $y$ is obtained from the order ideal corresponding to 
$x$ by adding a single element, labeled by $i$, say.  Then we set  
$\lambda(x,y) = i$.  This gives us a \snelling\ for $L=J(Q)$.
Figure~\ref{fig:distexample} shows a labeled poset and its lattice 
of order ideals with the appropriate edge-labeling. 
\begin{figure}
\center
\epsfxsize=40mm
\epsfbox{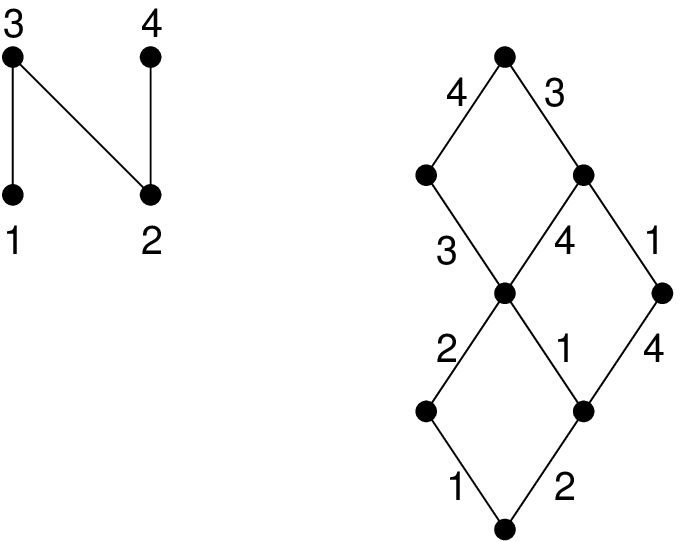}
\caption{}
\label{fig:distexample}
\end{figure}
\end{example}

\begin{example}\label{elshellable}
The posets shown in Figure \ref{fig:notsnellable} are seen to be 
EL-shellable.
However, it can be shown that neither of them is \snellable.  
Notice that the second poset, unlike the first, is a lattice.
It appears, together with this EL-labeling, in \cite{Si1}.   
\begin{figure}
\center
\epsfxsize=80mm
\epsfbox{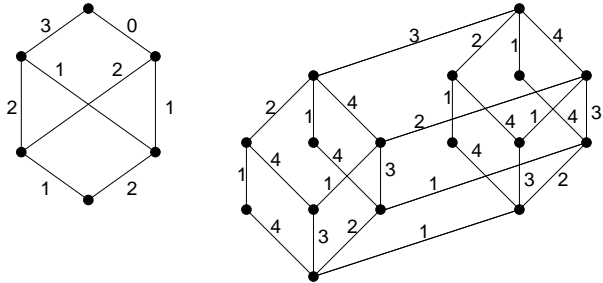}
\caption{Two posets that are EL-shellable but not \snellable}
\label{fig:notsnellable}
\end{figure}
\end{example}

\begin{example}\label{ss}
The set of supersolvable lattices is our final example and is also the 
example most relevant to the remainder of the paper.  
The following definition first appeared in \cite{St2}. 
\end{example}

\begin{definition}
A finite lattice $L$ is said to be \emph{supersolvable} if it contains
a maximal chain, called an \emph{M-chain} of $L$,  which together with any other chain in $L$ generates 
a distributive sublattice.
\end{definition}

We can label each such distributive sublattice by the method described in
Example \ref{distributive} in such a way that the M-chain receives the
increasing label 
$(1,2,\ldots,n)$.  As shown in \cite{St2}, this will assign 
a unique label to each edge of $L$ and the resulting global labeling
of $L$ is a \snelling.

Examples of supersolvable lattices include distributive lattices, 
the lattice $\Pi_n$ of partitions
of $[n]$, the lattice $NC_n$ of non-crossing partitions of $[n]$ and the
lattice $L(G)$ of subgroups of a supersolvable group $G$ (hence the
terminology).  The supersolvability of $\Pi_n$ and $L(G)$ was shown in
\cite{St2} while \cite{He} contains a proof that $NC_n$ is supersolvable. 

We are now in a position to state our first main result.  

\begin{theorem}\label{ss=snellable}
A finite graded lattice of rank $n$ is supersolvable if and only if it
is $\sn$ EL-shellable.
\end{theorem}

We will prove Theorem \ref{ss=snellable} in Section \ref{sec-proof1}.

\section{$\mathcal{H}_{\lowercase{n}}(0)$ actions}\label{sec-action}

Let $P$ be a finite graded poset of rank $n$ with $\hat{0}$ and $\hat{1}$.
Suppose $P$ has a \snelling\ $\lambda$.  Then to any maximal chain 
$\m : \hat{0}=x_0 < x_1 < \cdots < x_n=\hat{1}$ 
we can associate the permutation $\wm$ given by 
\[
\wm = 
(\lambda(x_0, x_1), \lambda(x_1, x_2), \ldots , \lambda(x_{n-1}, x_n) ).
\]
It is now natural to define the \emph{descent set of $\m$} to be 
the descent set of $\wm$ and the number of \emph{inversions} of $\m$ to
be the number of inversions of $\wm$. 
Suppose $\m$ has a descent at $i$.  By the snellability of $P$, 
there exists exactly one chain 
$\m' : \hat{0}=x_0 < x_1 < \cdots <x_{i-1} < x_{i}^{\prime} < x_{i+1}
< \cdots < x_n=\hat{1}$ differing only from $\m$ at rank $i$ and having
no descent at $i$.  
This suggests the following definition of functions $U_i:\mcp\to\mcp$.

\begin{definition}
Let $P$ be a finite graded poset of rank $n$ with $\hat{0}$ and $\hat{1}$
and with an $\sn$ EL-labeling.
Let $\m$ be a maximal chain of $P$.
We define $U_1,U_2,\ldots U_{n-1}:\mcp\to\mcp$ 
by $U_i(\m) = \m'$, where $\m'$ is the unique maximal chain
of $P$ differing only from $\m$ at possibly rank $i$ and having no 
descent at $i$.
\end{definition}

Under this definition, we see that the descent set of a maximal chain
$\m$ of $P$ can also be defined to be the set 
\begin{equation}\label{eq:descentset}
\{ i \in [n-1] : U_i(\m) \neq \m \}.
\end{equation}
This definition will be used later for posets $P$ where no \snelling\ is
defined. 

\begin{figure}
\center
\epsfxsize=12mm
\epsfbox{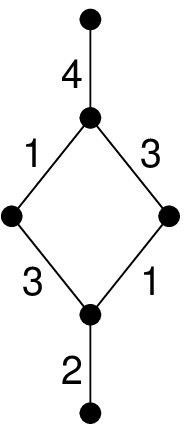}
\caption{}
\label{fig:parallelogram}
\end{figure}

Observe that $\omega_{\m'}$ is the same as $\wm$ except that the $i$th and
$(i+1)$st elements have been switched.  In other words, 
$\omega_{\m'} = \omega_\m s_i$ .
Figure \ref{fig:parallelogram}
shows an example for the case $n=4$.  Let $\m$ be the maximal chain to the
left.  It has a descent at $2$ and therefore $\m' = U_2(\m) \neq \m$.
The labels of $\m'$ are forced by the fact that $\m'$ does not have a descent 
at $2$.  We have that $\omega_{\m'} = \wm s_2$.

We see that the action of $U_1, U_2, \ldots, U_{n-1}$ has the following 
properties:
\begin{enumerate}
\item  It is a \emph{local} action, i.e., $U_i(\m)$  
agrees with $\m$ except possibly at
the $i$th rank.  Local actions on the maximal chains of a poset have been 
studied, for example, in \cite{He}, \cite{SS}, \cite{St3} and \cite{St4}. 
\item  ${U_i}^2 = U_i$ for $i=1,2,\ldots ,n-1$.  This differs from most of 
the local actions in the aforementioned papers which were
symmetric group actions and so satisfied ${U_i}^2 =  1$.
\item $U_i U_j = U_j U_i$ if $|i-j| \geq 2$.
\item $U_i U_{i+1} U_i = U_{i+1} U_i U_{i+1}$ for $i=1,2, \ldots ,n-2$.  
This requires the \snellable\ property and is left as an 
exercise for the reader. 
\end{enumerate}

Now we compare this to the definition of the 0-Hecke algebra  
$\mathcal{H}_n(0)$ as discussed in \cite{DHT}, \cite{DKLT} and \cite{KT}.

\begin{definition}
The 0-Hecke algebra $\mathcal{H}_n(0)$ of type $A_{n-1}$ 
is the $\mathbb{C}$-algebra generated
by $T_1, T_2, \ldots ,T_{n-1}$ with relations:
\begin{itemize}
\item[(i)]  ${T_i}^2 = -T_i$ for $i=1,2,\ldots ,n-1$.   
\item[(ii)] $T_i T_j = T_j T_i$ if $|i-j| \geq 2$.
\item[(iii)] $T_i T_{i+1} T_i = T_{i+1} T_i T_{i+1}$ for $i=1,2, \ldots ,n-2$. 
\end{itemize}
\end{definition}

We can extend the action of $U_1, U_2, \ldots, U_{n-1}$ on $\mcp$
to a linear action on $\mathbb{C}\mcp$, the complex vector space
with basis $\mcp$.
If we set \mbox{$T_i=-U_i$} then $U_1, U_2, \ldots, U_{n-1}$ generate the same 
$\mathbb{C}$-algebra and so we can now refer to our action on 
the maximal chains of $P$ as a \emph{local $\mathcal{H}_n(0)$ action.}
In 
\cite[{\S}3.9]{DHT}, Duchamps, Hivert and Thibon 
describe the special case of this action on distributive lattices. They
work in the language of linear extensions of a poset $Q$ which, as we
have seen, correspond to \snelling s of $J(Q)$.

Our action has one further very desirable property which we now discuss. 

\section{Good $\mathcal{H}_{\lowercase{n}}(0)$ actions}\label{sec-goodaction}

Before stating the fifth property, we must give some background, much of 
which is taken from the introduction in \cite{SS}.

Let $P$ be any finite graded poset of rank $n$ with $\hat{0}$ and $\hat{1}$
and let 
$S \subseteq [n-1]$. We let $\alpha_P(S)$ denote the number of 
chains in $P$ whose elements, other than $\hat{0}$ and $\hat{1}$,
have rank set equal to $S$.  In other words,
\[
\alpha_P(S) = \#\left\{\hat{0}<t_1<\cdots< t_{|S|} <\hat{1} :
		\left\{\rk(t_1),\ldots,\rk(t_{|S|}) \right\} = S
				\right\}.
\]
The function $\alpha_P:2^{[n-1]} \to \zz$ is called the 
\emph{flag f-vector} of $P$.  It contains equivalent information to
that of the \emph{flag h-vector} $\beta_P$ whose values are given by
\[
\beta_P(S) = \sum_{T \subseteq S} (-1)^{|S-T|}\alpha_P(T) .
\]

Ehrenborg in \cite[Def. 3]{Eh} suggested looking at the formal power series
(in the variables $x = (x_1, x_2, \ldots))$
\[
	F_P(x) = \sum_{\hat{0}=t_0 \leq t_1 \leq \cdots \leq t_{k-1}
			< t_k = \hat{1}}
		x_1^{\rk(t_0,t_1)} 
		x_2^{\rk(t_1,t_2)} 
		\cdots x_k^{\rk(t_{k-1},t_k)},
\]
where the sum is over all multichains from $\hat{0}$ to $\hat{1}$
such that $\hat{1}$ occurs exactly once.  
It is easy to see that the series $F_P(x)$ is homogeneous of degree $n$
and that it is a \emph{quasisymmetric function}, that is, for every
sequence $n_1, n_2, \ldots, n_m$ of exponents, the monomials 
$x_{i_1}^{n_1} x_{i_2}^{n_2} \cdots x_{i_m}^{n_m}$ and 
$x_{j_1}^{n_1} x_{j_2}^{n_2} \cdots x_{j_m}^{n_m}$  appear with 
equal coefficients whenever $i_1 < i_2 < \cdots <i_m$ and
$j_1 < j_2 < \cdots <j_m$.  
The series $F_P(x)$ can also be rewritten as 
\begin{equation}
F_P(x) = \sum_{S \subseteq [n-1]} \beta_P(S) L_{S,n}(x),
\label{eq:rewritefp}
\end{equation}
where $L_{S,n}(x)$ denotes Gessel's fundamental quasisymmetric function
\[
L_{S,n}(x) = \sum_{\genfrac{}{}{0pt}{}
		{1 \leq i_1 \leq i_2 \leq \cdots \leq i_n} 
		{i_j < i_{j+1} \mathrm{\: if\:} j \in S }}
	x_{i_1} x_{i_2} \cdots x_{i_n} , 
\]
which constitute a basis for the space of quasisymmetric functions of
degree $n$.  
The case when $F_P$ is 
a symmetric function is considered in \cite{SS} and \cite{St3}
and we wish, in a sense, to extend this to the case when $F_P$ is
a quasisymmetric function.
In our brief references to the symmetric function case,
we follow the notation of \cite{Ma}. 
The usual involution $\omega$ on symmetric functions given by
$\: \omega(s_{\lambda}) = s_{\lambda'} \:$ can be extended to the
ring of quasisymmetric functions by the definition \
$\omega(L_{S,n}) = L_{[n-1]-S\:,n}$. \
As in \cite[Exer. 7.94]{ec2}, 
where this extended definition appears, 
we leave it as an exercise to check that it restricts to
the ring of symmetric functions to give the usual $\omega$.

We now introduce some representation theory related to our
local $\mathcal{H}_n(0)$ action.
In the symmetric function case, certain classes of posets $P$ have been
found to have the property that 
\[
F_P(x) = \ch(\psi) \mbox{\ \ \ or\ \ \ } \omega F_P(x) = \ch(\psi)
\]
where $\psi$ denotes the character of some local symmetric group action 
and where $\ch(\psi)$ denotes its Frobenius characteristic 
as defined in \cite[{\S}I.7]{Ma}.
In extending these concepts to the $\mathcal{H}_n(0)$ case, we 
follow the definitions in \cite{DKLT} and \cite{KT}.
The representation theory of $\mathcal{H}_n(0)$ is studied by Norton
in \cite{No}.
There are known to be $2^{n-1}$ irreducible representations, all of 
dimension 1.  Since ${T_i}^2 = -T_i$, the irreducible representations 
are obtained by sending a set of generators to $-1$ and its complement 
to 0.  We will label these representations by subsets $S$ of $[n-1]$, 
and then the irreducible representation $\psi_S$ of $\mathcal{H}_n(0)$ is
defined by
\[
\psi_S(T_i) = \left\{ \begin{array}{rl}
			-1 & 	\mbox{if $i \in S$}, 	\\
			0 &	\mbox{if $i \not\in S$}.
		  \end{array}
		\right. 
\]
Therefore,
\[
\psi_S(U_i) = \left\{ \begin{array}{rl}
			1 & 	\mbox{if $i \in S$}, 	\\
			0 &	\mbox{if $i \not\in S$}.
		  \end{array}
		\right. 
\]
Hence the character of $\psi_S$, denoted by $\chi_S$, is given by
\[
\chi_S(U_{i_1}U_{i_2}\cdots U_{i_k})
	 = \left\{ \begin{array}{rl}
			1 & 	\mbox{if $i_j \in S$ for $j=1,2,\ldots,k$}\ , 	\\
			0 &	\mbox{otherwise}.
		  \end{array}
		\right. 
\]
We define its 
\emph{characteristic} by 
\[
\ch(\chi_S)=L_{S,n}(x),
\]
and we extend it to the set of all characters of representations of 
$\mathcal{H}_n(0)$ by linearity.
We let $\chi_P$ denote the character of the defining representation of
our local $\mathcal{H}_n(0)$ action on the space $\mathbb{C}\mcp$.

\begin{proposition} \label{property5}
Let $P$ be a finite \snellable\ graded 
poset of rank $n$ with $\hat{0}$ and $\hat{1}$.
Then the local $\mathcal{H}_n(0)$ action on the maximal chains of $P$ has
the property that 
\begin{equation} 
\omega F_P(x) = \ch(\chi_P).
\label{eq:fp=chi}
\end{equation} 
\end{proposition}

\begin{proof}
It is sufficient to show that the coefficient of $L_{S,n}$ for any 
$S \subseteq [n-1]$ is the same for both sides of \eqref{eq:fp=chi}. 
By \eqref{eq:rewritefp}, 
\[
\left[ L_{S,n} \right] \omega F_P(x) = \beta_P(S^c)
\]
where $S^c$ denotes $[n-1]-S$.

Let $J \subseteq [n-1]$ and let $\{i_1, i_2, \ldots, i_k\}$ be a multiset on
$J$ where each element of $J$ appears at least once.  Let $\m \in \mcp$.
If $U_i(\m) \neq \m$ for some $i \in [n-1]$ then $U_i(\m)$ has one less inversion 
than $\m$.  It follows that
$U_{i_1}U_{i_2}\cdots U_{i_k}(\m)=\m$ if and only if the descent set of $\m$
is disjoint from $J$.  Therefore 
\begin{eqnarray*}
\chi_P(U_{i_1}U_{i_2}\cdots U_{i_k}) & = & \#\left\{\m \in \mcp : \m  
					\mbox{ has no descents in $J$} \right\}
					\\
 & = & \sum_{S \supseteq J} \#\left\{\m \in \mcp : \m \mbox{ has descent set }S^c
					\right\} \\
 & = & \sum_{S \supseteq J} \beta_P(S^c) \mbox{\ \ \ by \cite[Thm. 2.2]{BGS} } \\
 & = & \sum_{S \subseteq [n-1]} \beta_P(S^c) 
		\chi_S(U_{i_1}U_{i_2}\cdots U_{i_k}).
\end{eqnarray*}
Thus
\[ \left[ L_{S,n} \right] \ch(\chi_P) 
= \left[ L_{S,n} \right] \ch \left( \sum_{S \subseteq [n-1]} \beta_P(S^c) 
					\chi_S \right)
= \beta_P(S^c) \]
as required.
\end{proof}

To summarize, we have that if 
$P$ is a finite \snellable\ graded poset 
of rank $n$, with $\hat{0}$ and $\hat{1}$,
then $P$ has a local $\mathcal{H}_n(0)$ action with the property that
$\omega F_P(x) = \ch(\chi_P)$.  Following \cite{St3}, we call such an action 
a \emph{good} $\mathcal{H}_n(0)$ action.
It is natural to ask what other types of posets have \ghna s.

\begin{figure}
\center
\epsfxsize=20mm
\epsfbox{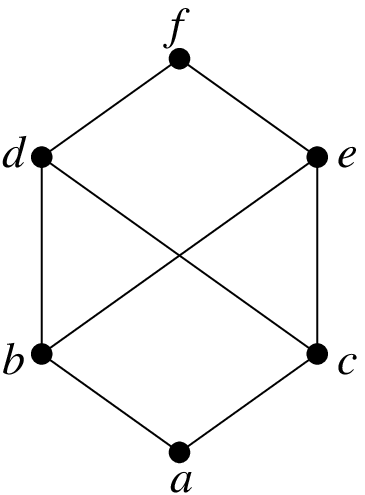}
\caption{}
\label{fig:goodaction}
\end{figure}

\begin{example}\label{goodaction}
Consider the poset $P$ shown in Figure \ref{fig:goodaction}.  
As stated in Example \ref{elshellable}, this poset is not \snellable.  
However, it does have a \ghna\ as described in the following table.
\begin{center}
\begin{tabular}{ccc}
$\m$ & $U_1(\m)$ & $U_2(\m)$ \\ \hline
$\m_1: a<b<d<f$ & $\m_3$ & $\m_2$ \\
$\m_2: a<b<e<f$ & $\m_4$ & $\m_2$ \\
$\m_3: a<c<d<f$ & $\m_3$ & $\m_4$ \\
$\m_4: a<c<e<f$ & $\m_4$ & $\m_4$ 
\end{tabular}
\end{center}

It is easy to check that this gives a local $\mathcal{H}_3(0)$ action.
We also have that 
\[
\omega F_P(x) = L_{\emptyset,3}+L_{\{1\},3}+L_{\{2\},3}+L_{\{1,2\},3} =
\ch(\chi_P).
\]
Therefore,  this poset has a \ghna.
\end{example}

\begin{definition}
A graded poset $P$ is said to be \emph{\shape-free} if it does not 
contain distinct elements $a, b, c$ and $d$ 
such that $a$ covers both $c$ and $d$, and such that $b$ covers both $c$ 
and $d$.
\end{definition}

In Section 3, we will prove our second main result:

\begin{theorem}\label{snellable=ghna}
Let $P$ be a finite graded \shape-free poset of rank $n$ with 
$\hat{0}$ and $\hat{1}$.
Then $P$ is $\sn$ EL-shellable if and only if $P$ has a \ghna.
\end{theorem} 

In particular, since lattices are \shape-free, we get the following 
immediate corollary.

\begin{corollary}\label{cor}
Let L be a finite graded lattice of rank $n$.  Then the following 
are equivalent:
\begin{enumerate}
\item	L is supersolvable,
\item	L is $\sn$ EL-shellable,
\item	L has a \ghna.
\end{enumerate}
\end{corollary}

\section{Snellable implies supersolvable} \label{sec-proof1}
Our main aim for this Section is to prove Theorem \ref{ss=snellable}.

Let $L$ be a finite graded lattice of rank $n$.  We showed in 
Example~\ref{ss} that if $L$ is supersolvable, 
then $L$ is \snellable.  Now 
we suppose that $L$ is \snellable\ and we wish to prove that $L$ is 
supersolvable.  We let $\m_0$ denote the unique maximal chain of $L$ 
labeled by the identity permutation.  Taking $\m_0$ to be our candidate 
M-chain, we let 
$L_{\chain}$ denote the sublattice of $L$ generated by $\m_0$ and any other
chain $\chain$ of $L$.  

It is shown in \cite[p.12]{Bi} and is easy to see that 
any sublattice of a distributive lattice is 
distributive.  If $\chain$ is a chain in $L$ that isn't maximal, then we can
extend it to a maximal chain $\m$ in at least one way.  Then $L_{\chain}$
is a sublattice of $L_\m$.  Therefore, it suffices to show that
$L_\m$ is distributive for all maximal chains $\m$.  Our approach will 
be to define two new posets, $Q_\m$ and $\jpm$, and
to show that 
\[
L_\m = Q_\m \cong \jpm. 
\]

We have seen that if $U_i(\m)$ differs from $\m$, then $U_i(\m)$ has
one less inversion than $\m$ and that $\omega_{U_i(\m)} = \omega_\m s_i$.
It follows that if $\m$ has $r$ inversions
then we can find a sequence
$U_{i_1}, U_{i_2}, \ldots, U_{i_r}$  
such that 
$U_{i_1}U_{i_2}\cdots U_{i_r}(\m)=\m_0$. 
We define $\Mm$, a subset of $\M(L)$, as follows:
\[
\Mm = \left\{ \m' \in \M(L) : \exists\ i_1, i_2, \ldots,i_r
\mbox{\ such that\ } m' = U_{i_1}U_{i_2}\cdots U_{i_r}(m) \right\}
\]
We label the elements of $\Mm$ as they are labeled in $L$.
We define $Q_\m$ to be the subposet of $L$ with elements
\[
\left\{ u \in L : u\in \m' \mbox{\ for some\ }\m' \in \Mm\right\}
\]
and with a partial order inherited from $L$.
$Q_\m$ can be thought of as the closure of $\m$ in $L$ under the operations
$U_1, U_2,\ldots,U_{n-1}$.  
We should note that it is not obvious that every maximal chain of 
$Q_\m$ is in $\Mm$.
We wish to obtain a clear picture of the
structure of $Q_\m$.

We are now ready to start the proof proper of Theorem \ref{ss=snellable}.
We break up the argument into a series of small steps.

\begin{step}\label{step:1}
Let $\m'$ and $\m''$ be distinct elements of
$\Mm$.  Then $\omega_{\m'} \neq \omega_{\m''}.$
\end{step}

Suppose that $\omega_{\m'} = \omega_{\m''}$.  Let 
$U_{i_1}, U_{i_2}, \ldots, U_{i_l}$  and 
$U_{j_1}, U_{j_2}, \ldots, U_{j_l}$ be 
sequences of minimal length such that
$\m' = U_{i_1}U_{i_2}\cdots U_{i_l}(\m)$ and 
$\m'' = U_{j_1}U_{j_2}\cdots U_{j_l}(\m)$.  
Then  
$s_{i_1}s_{i_2} \cdots s_{i_l}$ and 
$s_{j_1}s_{j_2} \cdots s_{j_l}$
are both reduced expressions for 
${\omega_{\m'}}^{-1}\omega_\m.$
By Tits' Word Theorem,
$s_{i_1}s_{i_2} \cdots s_{i_l}$ can thus be obtained from  
$s_{j_1}s_{j_2} \cdots s_{j_l}$ by a sequence of braid moves (i.e. replace
$s_is_{i+1}s_i$ by $s_{i+1}s_is_{i+1}$ or vice versa or replace
$s_is_j$ by $s_js_i$ if $|i-j|\geq2$.)
But by Properties 3 and 4 of the $U_i$ action,
$U_{i_1}U_{i_2}\cdots U_{i_l}(\m)$   
is invariant under braid moves.  We conclude that $\m'=\m''$, which is
a contradiction.  Therefore,
$\omega_{\m'} \neq \omega_{\m''}$.  
\subqed 

\begin{step}
Let $u \in Q_\m$.  Then there is a unique chain $\m_u \in \Mm$ that
has increasing labels between $\hat{0}$ and $u$ and between $u$ and
$\hat{1}$.
\end{step}

Choose any $\m' \in \Mm$ such that $u \in \m'$.  Suppose $u$ has rank $i$
in $L$.  Apply 
$U_1, U_2,\ldots,U_{i-1},U_{i+1},\ldots,U_{n-1}$ repeatedly to $\m'$ to
obtain $\m_u$.  The chain $\m_u$ is unique in $\Mm$ because it is 
unique in $L$. 
\subqed

\begin{step} \label{step:3}
To each point $u$ of $Q_\m$ we can associate the subset $\Lambda_u$ 
of $[n]$
consisting of the labels on any maximal chain of $[\hat{0}, u]$ in $L$.  Then 
any two distinct points of $Q_\m$ correspond to distinct subsets of $[n]$.
\end{step}

Let $u, v$ be distinct elements of $Q_\m$ and suppose 
$\Lambda_u = \Lambda_v$.  Then  
$\omega_{\m_u} = \omega_{\m_v}$,
contradicting Step \ref{step:1}.
\subqed

An important tool for the remainder of the proof will be the weak
order on permutations of $[n]$.

\begin{definition}
Let $v, w$ be permutations of $[n]$.  We say that $v \leq_R w$ if there exist 
$i_1, i_2,\dots,i_r$ such that 
$v=ws_{i_r}s_{i_{r-1}}\cdots s_{i_1}$ and
$ws_{i_r}\cdots s_{i_{k+1}}s_{i_k}$ has one less inversion than
$ws_{i_r}\cdots s_{i_{k+1}}$ for $k=1,2,\ldots,r$.
\end{definition}

It is known (see, for example, \cite[Prop. 2.5]{BW}) that $v \leq_R w$
if and only if ${INV}(v) \subseteq {INV}(w)$, where we define the
set of inversions of $v$, ${INV}(v)$, by
\[
{INV}(v) = \left\{ (v(j),v(i)) \in [n]\times [n] : i<j, v(i)>v(j) \right\}\ .
\]

\begin{step}
The labels on the elements of $\Mm$ consist of all those
permutations $\omega$ satisfying $\omega \leq_R \omega_\m$, each occurring
exactly once. 
\end{step}

Compare the definitions of $\Mm$ and $\leq_R$.  We see that if $\m' \in \Mm$
then $\omega_{\m'} \leq_R \omega_\m$ and if 
$\omega \leq_R \omega_\m$ then there exists $\m' \in \Mm$ satisfying 
$\omega_{\m'} = \omega$.  The fact that $\omega$ occurs only 
once follows from Step \ref{step:1}.
\subqed

\begin{step} \label{step:5}
Let $u, v \in Q_\m$.  We know that if $u \leq v$ then 
$\Lambda_u \subseteq \Lambda_v$.  Suppose $\Lambda_u \subseteq \Lambda_v$
for some elements $u, v$ of $Q_\m$.  Then $u \leq v$. 
\end{step}

Construct a permutation $\omega$ as follows:
\begin{itemize}
\item Let $\omega(1),\omega(2),\ldots,\omega(|\Lambda_u|)$ be the 
elements of $\Lambda_u$ taken in increasing order.
\item Let $\omega(|\Lambda_u|+1),\ldots,\omega(|\Lambda_v|)$ be the 
elements of $\Lambda_v-\Lambda_u$ taken in increasing order.
\item Let $\omega(|\Lambda_v|+1),\ldots,\omega(n)$ be the 
elements of $[n]-\Lambda_v$ taken in increasing order.
\end{itemize}
Then, since $u,v \in Q_\m$, we have that ${INV}(\omega) \subseteq 
{INV}(\omega_\m)$ and so $\omega \leq_R \omega_\m$.
Let $\m_{u,v}$ be the element of $\Mm$ satisfying 
$\omega_{\m_{u,v}} = \omega$.  By Step \ref{step:3}, $u$ and $v$ are both
elements of $\m_{u,v}$.  We conclude that $u \leq v$ in $Q_\m$.
\subqed

We can now exhibit a poset $P_{\omega_\m}$ such that 
$Q_\m \cong \jpm$.
Construct $P_{\omega_\m}$, a poset on $[n]$ with 
relation $\leq$ defined by 
$i \leq j$ if and only if 
$(i,j) \not\in {INV}(\omega_\m)$.  
For example, if $\omega_\m=2413$
we get the poset on the left in 
Figure \ref{fig:distexample}.  

\begin{step}
The map $\phi: Q_\m \to \jpm\ $ defined by $\phi(u)=\Lambda_u$ 
is an isomorphism.
\end{step}

Suppose $\Lambda_u$ has size $k$.
\begin{eqnarray*}
u \in Q_\m  & \Leftrightarrow &
	\Lambda_u=\left\{ \omega(1), \omega(2), \dots, \omega(k) \right\}
	\mbox{\ for some\ } \omega \leq_R \omega_\m
	\\
& \Leftrightarrow & 
	\Lambda_u=\left\{ \omega(1), \omega(2), \dots, \omega(k) \right\}
	\mbox{\ for some\ }  \omega  \mbox{\ satisfying\ } \\
&	&	{INV}(\omega) \subseteq {INV}(\omega_\m) \\
& \Leftrightarrow & 
	\mbox{$\Lambda_u$ is an order ideal of $P_{\omega_\m}$} \\
& \Leftrightarrow & 
	\Lambda_u \in \jpm.
\end{eqnarray*}
Therefore, $\phi$ is a well-defined bijection.
If $u$ and $v$ are elements of $Q_\m$, by Step \ref{step:5},
\begin{equation} \label{eq:uvrelation}
u \leq v \mbox{\ in\ } Q_\m \Leftrightarrow 
\Lambda_u \subseteq \Lambda_v
\Leftrightarrow \Lambda_u \leq \Lambda_v \mbox{\ in\ } \jpm
\end{equation}
as required.
\subqed

It follows from this that $Q_\m$, up to isomorphism, depends only
on $\wm$ and not even on the underlying lattice $L$.

\begin{step}
$Q_\m$ is a sublattice of $L$.
\end{step}

Let $u,v \in Q_\m$ with corresponding subsets $\Lambda_u$ and 
$\Lambda_v$, respectively.
Let $u \vee_{L} v$ denote the join of $u$ and $v$ in $L$ and let 
$u \vee_{Q_\m} v$ denote the join of $u$ and $v$ in $Q_\m$, which we now 
know is a lattice.
In $L$ we have that  
\[ u \vee_{Q_\m} v \ \geq\ u \vee_L v \]
since $Q_\m$ is a subposet of $L$.  But by \eqref{eq:uvrelation},
\[ \rk(u \vee_{Q_\m} v) = 
\left| \Lambda_u \cup \Lambda_v \right| \leq \rk(u \vee_L v)  \]
since there are maximal chains of $[\hat{0}, u \vee_L v]$ going through 
$u$ and others going through $v$.  Thus,
\[ u \vee_{Q_\m} v \ =\ u \vee_L v .\]
Similarly, 
\[ u \wedge_{Q_\m} v \ =\ u \wedge_L v .\]
\subqed

We have shown that $Q_\m$ is a distributive sublattice of $L$.
Furthermore, $L_\m$ is a sublattice of $Q_\m$ since 
$L_\m$ is a sublattice of $L$ and 
$Q_\m$ contains $\m$ and $\m_0$. 
We conclude that $L_\m$ is also distributive and hence $L$ is 
supersolvable.  
\qed

The astute reader will notice that, while we have shown that $L$ is 
supersolvable and that $L_\m \subseteq Q_\m$, we have not fulfilled our
promise to show that $L_\m = Q_\m$.   However, this follows from the
following lemma.

\begin{lemma}
For each element $\m'$ of $\Mm$, we have
$Q_{\m'} = L_{\m'}$.
\end{lemma}

\begin{proof}
Let $\m'$ be an element of $\Mm$ such that $\omega_{\m'}$ has $l$
inversions.
The proof is by induction on $l$ with the result being trivially true
for $l=0$.  
Since we know that $L_{\m'} \subseteq Q_{\m'}
\subseteq Q_\m$, it suffices to restrict 
our attention to $Q_\m$. We will label the elements of $Q_\m$ by 
their corresponding subsets of $[n]$.  By \eqref{eq:uvrelation},
join and meet in $Q_\m$ are just set union and set intersection, respectively.

\begin{figure}
\center
\epsfxsize=28mm
\epsfbox{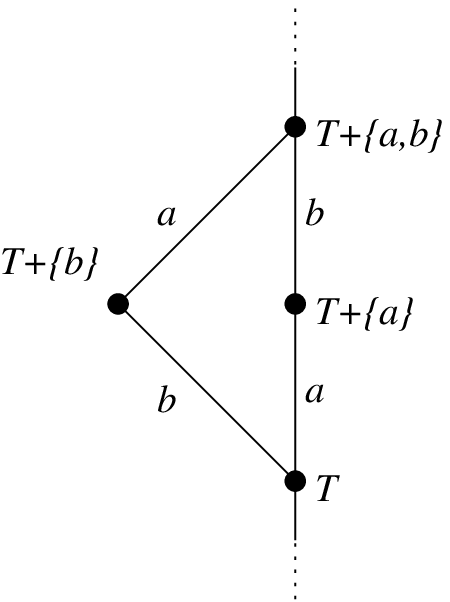}
\caption{}
\label{fig:sim}
\end{figure}
Referring to Figure \ref{fig:sim},
suppose $\m'$ is the vertical chain.
Suppose that $|T| = i-1$ and $a>b$ so that $\m'$ has a descent at rank $i$.
Now
\[ T + \{b\} = \left( \left(T + \{a,b\}\right) \cap \left( \{1,2,\ldots ,a-1\}
		\right) \right) \cup T \]
and $\{1,2,\ldots ,a-1\} \in \m_0$.  Therefore, $T+\{b\} \in L_{\m'}$ and so 
we get that $L_{U_i(\m')} \subseteq L_{\m'}$ as sets.
Suppose the descents of $\m'$ are at ranks $i_1, i_2, \ldots , i_k$.  
Then, as sets, 
\begin{eqnarray*}
Q_{\m'} & = & 
Q_{U_{i_1}(\m')} \cup 
Q_{U_{i_2}(\m')} \cup \cdots \cup
Q_{U_{i_k}(\m')} \cup \m' \\
& = & 
L_{U_{i_1}(\m')} \cup 
L_{U_{i_2}(\m')} \cup \cdots \cup
L_{U_{i_k}(\m')} \cup \m' \mbox{\ \ \ by induction}\\
& \subseteq & L_{\m'} .
\end{eqnarray*}
\end{proof}

\begin{example}
A \emph{non-crossing partition} of $[n]$ is a partition of $[n]$ into
blocks with the property that if some block $B$ contains $i$ and $k$ and
some block $B'$ contains $j$ and $l$ with $i<j<k<l$ then $B=B'$.
We order the set of non-crossing partitions by \emph{refinement}: if 
$\mu$ and $\nu$ are non-crossing partitions of $[n]$ we say
that $\mu \leq \nu$ if every block of $\mu$ is contained in 
some block of $\nu$.  The resulting poset $NC_n$, which is a subposet
of the 
lattice $\Pi_n$ of partitions of $[n]$, is itself a lattice
and has been studied extensively.  More information on $NC_n$ can be found 
in Simion's survey article \cite{Si2} and the references given there.

$\Pi_{n+1}$ was shown to be supersolvable in \cite{St2} and so can be
given a \snelling\ $\lambda$ as in Example~\ref{ss}.  We can choose the
M-chain to be the maximal chain consisting of the bottom element
and those partitions of $[n+1]$ whose only 
non-singleton block is $[i]$ where
$2 \leq i \leq n+1$.
In the literature, $\lambda$ is often seen in the following form, which
can be shown to be equivalent.
If $\nu$ covers $\mu$ in $\Pi_{n+1}$, then $\nu$ is obtained from $\mu$
by merging two blocks $B$ and $B'$ of $\mu$.
We set
\[
        \lambda(\mu, \nu) = \max \left\{ \min B, \min B'\right\} -1 .
\]
It was observed by A. Bj{\"{o}}rner and P. Edelman
in \cite{Bj}
that $\lambda$ restricts to $NC_{n+1}$ to give an EL-labeling for $NC_{n+1}$. 
In fact, it is readily checked that 
we get a \snelling\ for 
$NC_{n+1}$.  Theorem \ref{ss=snellable} now gives a new proof of the
supersolvability of $NC_{n+1}$.
\begin{figure}
\center
\epsfxsize=100mm
\epsfbox{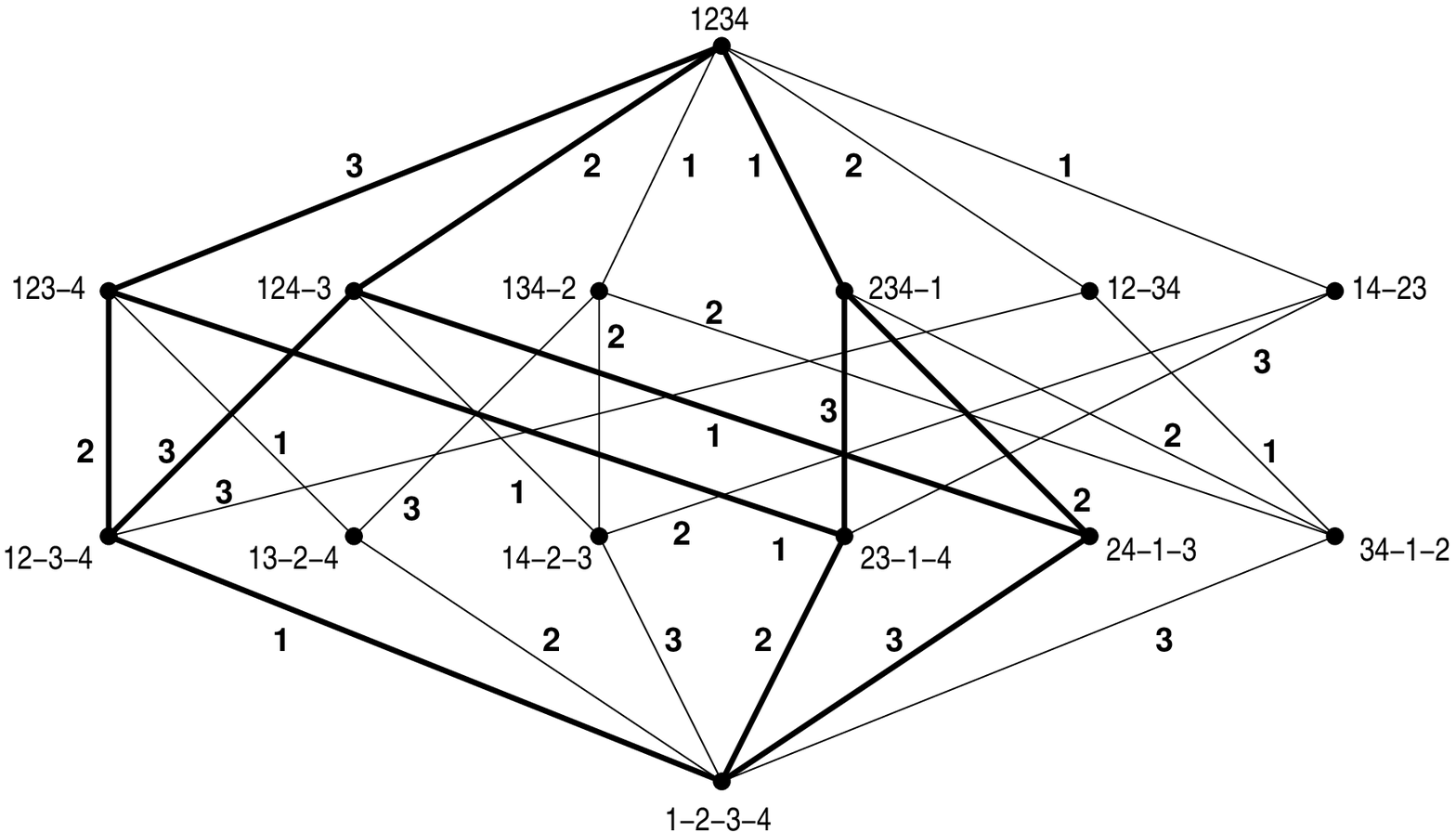}
\caption{$NC_4$ with \snelling}
\label{fig:nc}
\end{figure}
Figure~\ref{fig:nc} shows $NC_4$ with $L_\m=Q_\m$ highlighted 
for when $\m$ is the maximal chain $\hat{0} < \mbox{24-1-3} < 
\mbox{234-1} < \hat{1}$.
In this case, $P_{\omega_\m}$ is just 3 incomparable elements and so 
$\jpm = B_3 \cong Q_\m$.  

\end{example}

\section{Lattice Snellings and Good $\mathcal{H}_{\lowercase{n}}(0)$ Actions} 
\label{sec-proof2}
Our main aim for this Section is to prove Theorem \ref{snellable=ghna}.

Recall that $P$ denotes a finite graded \shape-free poset 
of rank $n$ with 
$\hat{0}$ and $\hat{1}$.
We suppose that $P$ has a \ghna\ and we let $\chi_P$ denote the character
of the defining representation of this action on the space $\mathbb{C}\mcp$.
In other words, we suppose that there 
exist functions 
$U_1, U_2, \ldots, U_{n-1}: \mcp \to \mcp$ satisfying the 
following properties:
\begin{enumerate}
\item The action of $U_1, U_2, \ldots , U_{n-1}$ is local.
\item ${U_i}^2 = U_i$ for $i=1,2,\ldots ,n-1$.  
\item $U_i U_j = U_j U_i$ if $|i-j| \geq 2$.
\item $U_i U_{i+1} U_i = U_{i+1} U_i U_{i+1}$ for $i=1,2, \ldots ,n-2$.  
\item $\omega F_P(x) = \ch(\chi_P)$.
\end{enumerate}
As we have previously suggested, 
given any maximal chain $\m$ of $P$, we define the descent set of $\m$ to
be the set 
\[
\{ i \in [n-1] : U_i(\m) \neq \m \}.
\]
We wish to show that $P$ is \snellable.
The following approach was suggested by R. Stanley.  
Suppose $P$ has a unique maximal chain $\m_0$ with empty descent set.
Given a maximal chain $\m$ of $P$,
suppose we can find 
$U_{i_1}, U_{i_2}, \ldots , U_{i_r}$ with $r$ minimal such that 
\begin{equation} \label{eq:noloops}
U_{i_1}U_{i_2}\cdots U_{i_r}(\m) = \m_0 .
\end{equation}
Then to $\m$ we associate the permutation $\wm = 
s_{i_1}s_{i_2}\cdots s_{i_r}$ 
and we label the edges of $\m$ 
by $\wm(1), \wm(2), \ldots, \wm(n)$ from bottom to top.
Our proof of the validity of this approach divides into four main 
parts.  
The first task is to show that $\m_0$ exists and is unique.
The next is to show that, given $\m$, we can always 
find 
$U_{i_1}, U_{i_2}, \ldots , U_{i_r}$ satisfying \eqref{eq:noloops}.
The third task is to show that $\wm$ is well-defined. 
Finally, we must show that this gives a \snelling\ for $P$.

\begin{definition}
Given maximal chains $\m$ and $\m'$ of $P$, we say that the expression
\mbox{$U_{i_1}U_{i_2}\cdots U_{i_r}(\m) = \m'$} is \emph{\restless} if
$U_{i_r}(\m) \neq \m$ and if 
\[
\mbox{$U_{i_j}U_{i_{j+1}}\cdots U_{i_r}(\m) \neq U_{i_{j+1}}\cdots U_{i_r}(\m)$
\ \ for $j=1,2,\ldots ,r-1$.} 
\]
\end{definition}
We say that two sequences 
$U_{i_1}U_{i_2}\cdots U_{i_r}$ and 
$U_{j_1}U_{j_2}\cdots U_{j_r}$ are in the same \emph{braid class} if we can get
from one to the other by applying Properties 3 and 4 repeatedly.  
It can be readily checked that if 
$U_{i_1}U_{i_2}\cdots U_{i_r}$ and 
$U_{j_1}U_{j_2}\cdots U_{j_r}$ are in the same braid class
and if 
\mbox{$U_{i_1}U_{i_2}\cdots U_{i_r}(\m) = \m'$} is \restless, then
\mbox{$U_{j_1}U_{j_2}\cdots U_{j_r}(\m) = \m'$} is \restless.
Here we use the \shape-free property of $P$.

To every sequence $i_1, i_2, \ldots ,i_r$ such that 
\mbox{$U_{i_1}U_{i_2}\cdots U_{i_r}(\m) = \m'$}, we can associate a counting 
vector of length $n-1$ where the $j$th coordinate equals the number  
of times that $i_j$ appears in the sequence $i_1, i_2, \ldots, i_r$.
We say that the expression
\mbox{$U_{i_1}U_{i_2}\cdots U_{i_r}(\m) = \m'$} is \emph{lexicographically
minimal} (or \emph{lex.\ minimal} for short) if no sequence 
$U_{j_1}U_{j_2}\cdots U_{j_r}$ in the braid class of
$U_{i_1}U_{i_2}\cdots U_{i_r}$ 
and satisfying
\mbox{$U_{j_1}U_{j_2}\cdots U_{j_r}(\m) = \m'$} has a lexicographically
less counting vector.  

The following result will help us to complete our first two tasks.

\begin{lemma}\label{noloops}
Let $\m'$ be any maximal chain of $P$.  Suppose
\mbox{$U_i(\m')\neq \m'$}.  
Then there do not exist $i_1, i_2, \ldots, i_r$
satisfying 
\mbox{$U_{i_1}U_{i_2}\cdots U_{i_r}U_i(\m')=\m'$.} 
\end{lemma}

\begin{proof}
Suppose there exist $i_1, i_2, \ldots, i_r$ 
satisfying 
\mbox{$U_{i_1}U_{i_2}\cdots U_{i_r}U_i(\m')=\m'$.} 
It suffices to consider the case when 
\mbox{$U_{i_1}U_{i_2}\cdots U_{i_r}U_i(\m')=\m'$} 
is \restless\ and lex.\ minimal.
Let $l \in [n-1]$ denote the minimum element of the sequence 
$i_1, i_2, \ldots,i_r, i$.  Since our equation is \restless\, $U_l$
must occur at least twice in the sequence.  

Take any pair of 
$U_l$ appearances with no $U_l$ between them.  If we had no $U_{l+1}$
between them, we could apply Property 3 until we had an appearance of
$U_lU_l$, contradicting the \restless\ property since ${U_l}^2 = U_l$. 
If there is just one $U_{l+1}$ between them, we can apply Property 3 to
get $U_lU_{l+1}U_l$ appearing and then apply Property 4 to get 
$U_{l+1}U_lU_{l+1}$, contradicting the lex.\ minimal property.
We conclude that, 
between the two appearances of $U_l$,  
there are at least two appearances of $U_{l+1}$.
Choose any two of these appearances of $U_{l+1}$
that don't have another $U_{l+1}$ between them and apply 
the same argument to show that there must be at least two appearances of 
$U_{l+2}$ between them.
Repeating this process, we eventually get $U_iU_i$ appearing, 
yielding a contradiction.
\end{proof}

More generally, we can apply the same argument to prove the following 
statement:

\begin{lemma}\label{lowestonce}
Suppose 
$U_{i_1}U_{i_2}\cdots U_{i_r}(\m) = \m_0$ is \restless\ and lex.\ minimal.
Let $l$ denote the minimum element of the sequence 
$i_1, i_2, \ldots,i_r$.  Then $U_l$ appears exactly once and for 
$l < i \leq n-1$, there must be an appearance of $U_{i-1}$ between any
two appearances of $U_i$.
\end{lemma}

The following result is essentially a rephrasing of Property 5 of
our \ghna\ into more amenable terms.

\begin{proposition}\label{rephrasing}
For all $S \subseteq [n-1]$, $\alpha_P(S)$ equals the number of 
maximal chains of $P$ with descent set contained in $S$.
\end{proposition}

\begin{proof}
We know that 
\[ \chi_P = \sum_{S\subseteq [n-1]}b_{P,S}\chi_S \]
for some set of coefficients $\{b_{P,S}\}_{S \subseteq [n-1]}$ and hence 
\[ \ch(\chi_P) = \sum_{S\subseteq [n-1]}b_{P,S}L_{S,n}.\]
By \eqref{eq:rewritefp} and Property 5, we see that $b_{P,S} = \beta_P(S^c)$.

Now let $J=\{i_1, i_2,\ldots,i_k\}\subseteq [n-1]$.  
Then
\begin{eqnarray*}
\sum_{S \supseteq J}\beta_P(S^c) & = &
\sum_{S \subseteq [n-1]}\beta_P(S^c)\chi_S(U_{i_1}U_{i_2}\cdots U_{i_k}) \\
& = & \chi_P(U_{i_1}U_{i_2}\cdots U_{i_k}) \\
& = & \#\left\{\m \in \mcp : \m  \mbox{ has no descents in $J$} \right\} 
\end{eqnarray*}
by Lemma \ref{noloops}.  Therefore,
\[
\sum_{S \supseteq J}\beta_P(S^c) 
 = \sum_{S \supseteq J} \#\left\{\m \in \mcp : \m \mbox{ has descent set }S^c
					\right\}. 
\]
Since this holds for all $J \subseteq [n-1]$, we get that
\[
\beta_P(S) = \#\left\{\m \in \mcp : \m \mbox{ has descent set }S \right\} 
\]
for all $S \subseteq [n-1]$.
By Inclusion-Exclusion,
this is equivalent to 
\[
\alpha_P(S) = \#\left\{\m \in \mcp : \m \mbox{ has descent set contained in }S \right\}. 
\]
\end{proof}

In particular, setting $S=\emptyset$, we see that $P$ has exactly one 
maximal chain, which we denote by $\m_0$, with no descents. 
Also, given a maximal chain $\m$ of $P$,
by Lemma~\ref{noloops} and the finiteness of $P$, 
we can find 
$U_{i_1}, U_{i_2}, \ldots , U_{i_r}$ with $r$ minimal such that 
$U_{i_1}U_{i_2}\cdots U_{i_r}(\m) = \m_0$ .
This completes our first two tasks.

Given any maximal chain $\m$ of $P$,
we consider the braid classes of the set of sequences 
$U_{i_1},U_{i_2},\ldots,U_{i_r}$ such that
$U_{i_1}U_{i_2}\cdots U_{i_r}(\m) = \m_0$ is \restless.
Our next task is to show that there is only one such braid class. 
Every braid class contains at least one element 
$U_{i_1},U_{i_2},\ldots,U_{i_r}$ such that
$U_{i_1}U_{i_2}\cdots U_{i_r}(\m)~=~\m_0$ is \restless\ and lex.\ minimal.
For such an element, the minimum, $l$, of $i_1,i_2,\ldots,i_r$ is the lowest 
rank for which $\m \neq \m_0$, by Lemma \ref{lowestonce}.
It follows that $l$ is the same for all the braid classes.  
It suffices to consider the case when $l=1$. 

The following result is central to our proof that there is just one 
braid class.

\begin{lemma} \label{endsini}
Suppose that the expressions
$U_{i_1}U_{i_2}\cdots U_{i_r}(\m) = \m_0$ and 
\linebreak
$U_{j_1}U_{j_2}\cdots U_{j_s}(\m)=\m_0$ are both \restless.
Then there exists an element of the braid class of 
$U_{i_1}U_{i_2}\cdots U_{i_r}$ and an element of the braid class of 
$U_{j_1}U_{j_2}\cdots U_{j_s}$ both ending on the right with the same
$U_i$.
\end{lemma}

\begin{proof}
Suppose 
$U_{i_1}U_{i_2}\cdots U_{i_r}$ and
$U_{j_1}U_{j_2}\cdots U_{j_s}$ are in different braid classes.  
Without loss of generality, we take them both to be lex. minimal.
If $U_1$ can be moved to the right-hand end in both by applying 
Property~3, then there's nothing to prove.
Suppose, by applying Property 3, that $U_1$ can be brought to the right 
end in one sequence but not in the other.  
Then $P$ must have the edges shown in Figure \ref{fig:proof2}, where 
$\m$ and $\m_0$ are the maximal chains on the left and right, respectively.
\begin{figure}
\center
\epsfxsize=18mm
\epsfbox{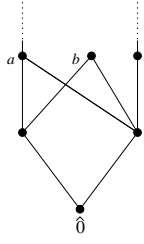}
\caption{A portion of $P$}
\label{fig:proof2}
\end{figure}
We see that we get a contradiction with the \shape-free property unless
$a=b$. In this case, $U_2$ appears at least twice in the latter sequence 
to the right of  
the unique appearance of $U_1$, contradicting Lemma 
\ref{lowestonce}.  We conclude that $U_1$ can't be brought to the
right end in either sequence.  Now we consider that portion of each
sequence to the right of the unique $U_1$.  By the same logic, the
maximal chains we get when we apply these portions to $\m$ must have
the same element at rank 2.

Consider the unique $U_2$ in each of these portions.  By a 
similar argument, we conclude that either we've nothing to prove or 
else $U_2$ can't be brought to the right of either sequence by applying
Property 3.  In the latter case, we consider the 
portion of each sequence to the 
right of the unique $U_2$. The maximal chains we get when we apply these
portions to $\m$ must have the same element at rank 3.  Repeating the same 
argument, we are eventually reduced to the case where $U_i$ is the 
element at the right end of both sequences, for some $i$.
\end{proof}

\begin{proposition}
If 
$U_{i_1}U_{i_2}\cdots U_{i_r}(\m) = \m_0$ and 
$U_{j_1}U_{j_2}\cdots U_{j_s}(\m) = \m_0$ are both \restless\, then 
$s_{i_1}s_{i_2}\cdots s_{i_r} = s_{j_1}s_{j_2}\cdots s_{j_s}$. 

\end{proposition}

\begin{proof}
It suffices to prove the result in the case when $r$ is as small as
possible.  We prove the result by induction on $r$, the result being
trivially true when $r=0$.  

For $r>0$, by the previous lemma, there exists an element 
$U_{I_1}U_{I_2}\cdots U_{I_{r-1}}U_i$ of the braid class of
$U_{i_1}U_{i_2}\cdots U_{i_r}$ and an element 
$U_{J_1}U_{J_2}\cdots U_{J_{s-1}}U_i$ of the braid class of
$U_{j_1}U_{j_2}\cdots U_{j_s}$.  Consider $U_i(\m)$.
By the induction hypothesis,
\[
s_{I_1}s_{I_2}\cdots s_{I_{r-1}} =
s_{J_1}s_{J_2}\cdots s_{J_{s-1}}.
\]
Therefore, since permutations are invariant under braid moves, 
\[
s_{i_1}s_{i_2}\cdots s_{i_r} =  
s_{I_1}s_{I_2}\cdots s_{I_{r-1}}s_i =
s_{J_1}s_{J_2}\cdots s_{J_{s-1}}s_i =
s_{j_1}s_{j_2}\cdots s_{j_s}.
\]
\end{proof}

Finally, we can make the following definition:

\begin{definition}
If 
$U_{i_1}U_{i_2}\cdots U_{i_r}(\m) = \m_0$ is \restless\ then we define
\[
\wm = s_{i_1}s_{i_2}\cdots s_{i_r}. 
\]
\end{definition}

For every maximal chain $\m$ of $P$, we
label the edges of $\m$ from bottom to top by 
$\wm(1), \wm(2), \ldots, \wm(n)$. 
Our final task is to show that this gives an edge-labeling, and in
particular a \snelling, for $P$.
We divide the proof into a number of small steps.

\setcounter{step}{0}
\begin{step}
If
$U_{i_1}U_{i_2}\cdots U_{i_r}(\m) = \m_0$ is \restless\ then 
$\wm = s_{i_1}s_{i_2}\cdots s_{i_r}$ is a reduced expression.  Furthermore,
if 
$\wm = s_{j_1}s_{j_2}\cdots s_{j_r}$ is another reduced expression,
then 
$U_{j_1}U_{j_2}\cdots U_{j_r}(\m) = \m_0$ is \restless.
\end{step}

The first assertion follows from the fact that if 
$\wm = s_{i_1}s_{i_2}\cdots s_{i_r}$ is not reduced then we can apply 
a sequence of braid moves to get $s_is_i$ appearing.  This contradicts
the \restless\ property. 
The second assertion follows from Tits' Word Theorem. 
\subqed

\begin{step} \label{step:2.2}
The permutation $\wm$ has a descent at $i$ if and only if 
$U_i(\m) \neq \m$.
In this case, $\omega_{U_i(\m)}$ is the same as $\wm$ except that
the $i$th and $(i+1)$st elements have been switched, removing the descent.
\end{step}

\begin{eqnarray*}
& U_i(\m) \neq \m \\
\Leftrightarrow & U_{i_1}U_{i_2}\cdots U_{i_r}U_i(\m)=\m_0 
\mbox{ is \restless\ for some $i_1, i_2, \ldots,i_r$}
\\
\Leftrightarrow & 
s_{i_1}s_{i_2}\cdots s_{i_r}s_i = \wm 
\mbox{ is a reduced expression for some $i_1, i_2, \ldots,i_r$} \\
\Leftrightarrow & \wm s_i \mbox{ has one less inversion than }  \wm\\
\Leftrightarrow & \wm \mbox{ has a descent at } i.
\end{eqnarray*}
When $U_i(\m)\neq \m$ and $\wm=s_{i_1}s_{i_2}\cdots s_{i_r}s_i$ 
is reduced we see
that 
$\omega_{U_i(\m)}=s_{i_1}s_{i_2}\cdots s_{i_r}$, yielding the 
second statement.
\subqed

\begin{step} \label{step:2.3}
Let $S \subseteq [n-1]$.  Then every chain in $P$ with rank set equal to  
$S$ has exactly one extension to a maximal chain of $P$ with 
descent set contained in $S$.
\end{step}

Given any chain $\chain$ with rank set $S$, let $\m$ be any extension
of $\chain$ to a maximal chain in 
$P$.  Apply $U_i$ for $i \not\in S$ repeatedly to $\m$.  By Step \ref{step:2.2}, 
this will eventually yield an extension of $\chain$
which is a maximal chain with descent set contained in $S$.
Therefore, every chain with rank set $S$ has at least one such extension.
We get
\[
\alpha_P(S)  \leq  
\#\left\{\m \in \mcp : \m \mbox{ has descent set contained in $S$} 
\right\}. 
\]
However, by Proposition \ref{rephrasing},
\[
\alpha_P(S) = \#\left\{\m \in \mcp : \m \mbox{ has descent set contained in $S$}
\right\}. 
\]
Thus $\chain$ has exactly one extension to a maximal chain of $P$ with
descent set contained in $S$.
\subqed

\begin{step}
For every maximal chain $\m$ of $P$, 
labeling the edges of $\m$ from bottom to top by 
$\wm(1), \wm(2), \ldots, \wm(n)$ gives an edge-labeling for $P$.
\end{step}

Let $x,y \in P$ be such that $y$ covers $x$ and let $\m$ and $\m'$
be maximal chains of $P$ containing both $x$ and $y$. 
Define $S = [\rk(x),\rk(y)]$ and let $\m_{(x,y)}$ denote the unique 
extension of $x<y$ to a maximal chain with descent set contained in $S$.
By applying $U_i$ for $i \not\in S$ repeatedly to $\m$, we can reach
$\m_{(x,y)}$. By Step \ref{step:2.2}, 
$\m$ and $\m_{(x,y)}$ give the same label to the edge $(x,y)$.
Similarly, 
$\m'$ and $\m_{(x,y)}$ give the same label to the edge $(x,y)$.
Therefore, 
$\m$ and $\m'$ give the same label to the edge $(x,y)$ and so we have an
edge labeling for $P$. 
\subqed

\begin{step}
This edge-labeling is a \snelling\ for $P$.
\end{step}

Let $x,y \in P$ be such that $x < y$.  
Let 
\[S = \left[n-1\right] - \left\{\rk(x)+1, \rk(x)+2,\ldots,\rk(y)-1\right\}\]
in Step \ref{step:2.3}.
The fact that the interval $[x,y]$ has exactly one increasing maximal chain
follows from Step \ref{step:2.3} and the fact that we now have an edge-labeling. 
Every maximal chain is labeled by a permutation by definition.
Therefore, $P$ is \snellable, proving Theorem \ref{snellable=ghna}.
\qed

\begin{remark}
Theorem \ref{snellable=ghna} does indeed contain information not contained in 
Corollary \ref{cor}, in that there exist finite graded 
\shape-free posets with $\hat{0}$ and $\hat{1}$ that are \snellable\ but
are not lattices.
For example, take the lattice $B_4$ with a snelling as described in 
Example \ref{bn}.  Now delete the edge (\{3,4\},\{2,3,4\}) in the Hasse
diagram of $B_4$ to form the Hasse diagram of a new poset.  
We can check that the new poset has the desired properties.  
\end{remark}

\begin{figure}
\center
\epsfxsize=35mm
\epsfbox{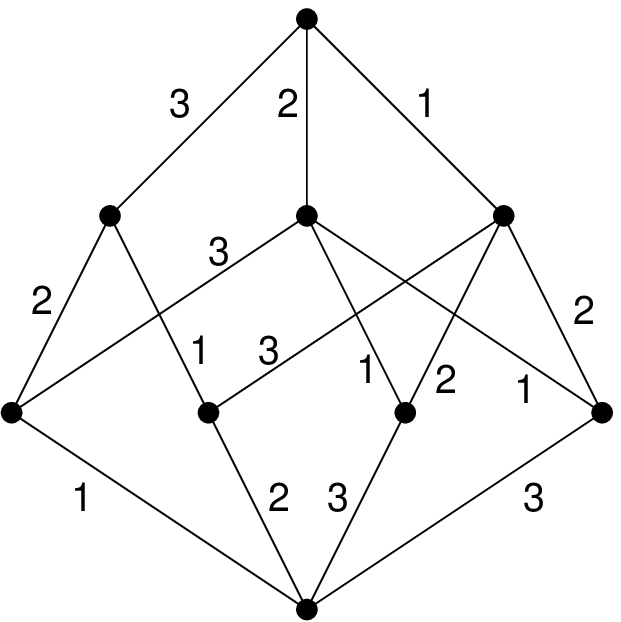}
\caption{}
\label{fig:counterex}
\end{figure}

It seems that we have fully answered the question of finite graded posets with
$\hat{0}$ and $\hat{1}$ in the \shape-free case.  What can we say about such
posets that are not \shape-free?  In Example \ref{goodaction} 
we saw a poset with
a \shape\ that has a \ghna\ but which is not \snellable.  On the 
other hand, Figure 
\ref{fig:counterex} shows a finite graded poset with $\hat{0}$ and $\hat{1}$
that has a \shape\ but which is still \snellable.  

This suggests the following question.

\begin{question}
Let $\mathcal{C}$ denote
the class of finite graded posets with $\hat{0}$, $\hat{1}$ 
and a \ghna.  
Is there some ``nice'' characterization of $\mathcal{C}$, possibly in 
terms of edge-labelings?
\end{question}

\section*{Acknowledgments}
The results of this paper all came from questions posed 
by Richard Stanley.  I am very grateful 
to him for this and for
many interesting and helpful discussions, especially regarding
\ghna s.
I am also grateful to Paul Edelman and Patricia Hersh 
for enthusiastic help with obtaining references and to Noam 
Elkies for technical help.

\end{document}